\documentclass[12pt]{article}
\usepackage[a4paper,nohead,left=20mm,right=20mm,top=20mm,bottom=20mm]{geometry}
\usepackage{local}

\begin{document}

\title{Classification of $k$-tangle projections using cascade representation}

\author{A. Bogdanov, V. Meshkov, A. Omelchenko, M. Petrov\\
St Petersburg Polytechnic University, St Petersburg, Russia}

\abstract{The paper addresses the $k$-tangle enumeration problem. We introduce a notion of
cascade diagram for $k$-tangle projections. An effective enumeration algorithm for projections is
proposed based on cascade representation. Tangles projections with up to 12 crossings are
tabulated. We provide also pictures of alternating $k$-tangles with 5 crossing or less.}

\maketitle

\section*{Introduction}

Tangles (see~\picref{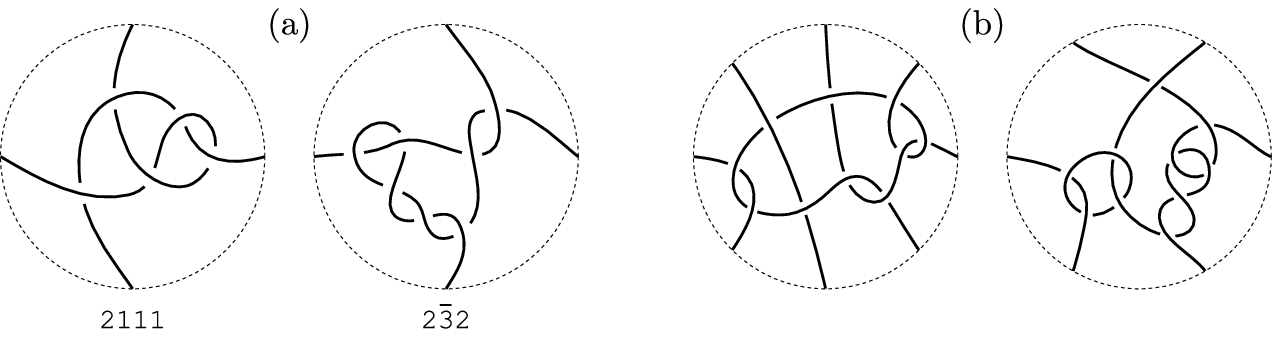}) were introduced by John Conway \cite{Conway1970} as an
instrument for the description and classification of knots and links. According to up-to-date
definitions~\cite{Cromwell2004}, Conway's tangles are 2-tangles. A natural generalization of
2-tangles are $k$-tangles or $k$-string tangles (\picref{conway.eps}b), where $k$ is one-half of
the number of ``legs''.

As in the case of knots and links we can set the classification problem for 
$k$-tangles. This problem is only solved for some specific subclasses of
2-tangles~\cite{Kanenobu2003,Kauffman2004} (e.\,g.~rational tangles) which are most-used for
studying knots. The general case of $k$-tangles is interesting in itself. Furthermore it is
possible to use $k$-tangles for representing knots and links in a ``thickened'' closed oriented
surface of arbitrary genus (virtual links~\cite{Kauffman1999,Kuperberg2003}). So $k$-tangle
classification can be a basis for classifying virtual links.


\centerpic{conway.eps}{2-tangles (Conway's notation) (a), $k$-tangles in general case (b)}

In this paper we introduce a notion of cascade diagram and propose a classification algorithm for
$k$-tangles based on the cascade representation. We represent cascade diagrams by a code that is
almost as compact as DT-code (Dowker--Thistlethwaite) for knots \cite{Dowker1983}. This cascade
code actually is a recipe how to draw the projection so there is no problem with verification of
realizability of a given code.

The first step in most of knots enumeration algorithms is enumeration of projections. Having the
full set of projections with $n$ crossings, it is easy to generate all the knots (links, tangles)
by arranging under- and over-crossing in each projection by all $2^n$ possible ways
\cite{Hoste2006}.

Projections enumeration is related directly to enumeration of alternating knots (links, tangles).
For any projection it is possible to arrange crossings in such a way that ``under'' and ``over''
will alternate when traveling along each component of the link (tangle). It is known that a
(reduced) alternating diagram of a given link has minimal number of crossings (Tait conjecture).
Any projection determines uniquely (up to mirror symmetry) an alternating knot (link, tangle).
However, the set of projections with $n$ crossings is wider than the set of alternating
knots (links, tangles): 
all alternating diagrams obtained from a given one by a series of flipes (\picref{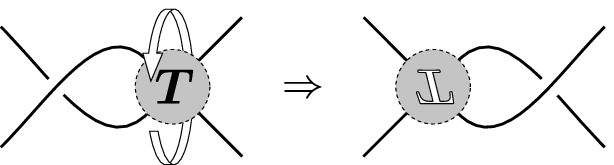})
correspond to the same knot (link, tangle).

\centerpic{flyp.eps}{A flip}

In the paper we will concentrate first of all on classification of tangle projections.


\section{Basic definitions}\label{defs}

\begin{dfn} A proper embedding of a disjoint union of $k$ arcs and (possibly) some
number of circles into the standard $3$-ball $B^3$ is called a {\it $k$-tangle.} The $2k$
endpoints are fixed at the points $(\cos\pi i/k,\sin\pi i/k,0)$, $i\in\{0,1,\dots,2k{-}1\}$ on
the boundary sphere~(\picref{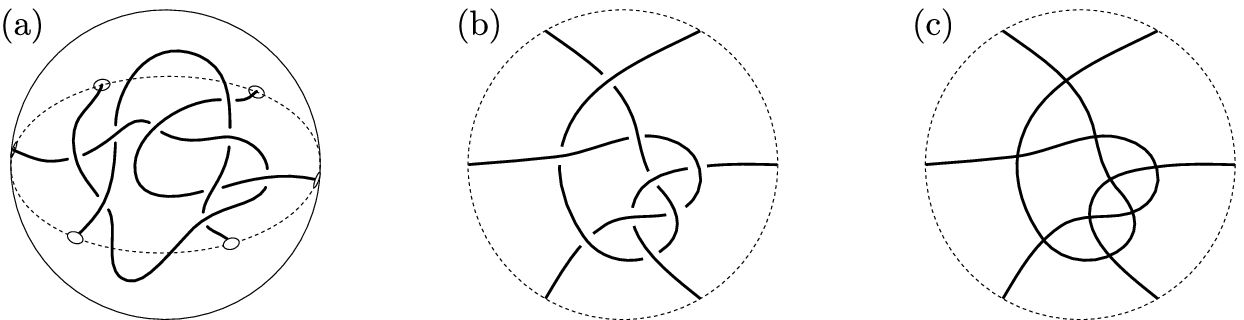}). $k$-tangle is considered up to the dihedral symmetry
group $D_{2k}$.\end{dfn}

\begin{dfn} Tangles $T_1$ and $T_2$ are {\it equivalent\/} if there is an isotopy i.\,e.~a
smooth family $f_t\,{:}\;B^3\to B^3$, $t\in[0,1]$ of smooth maps that keep the boundary sphere
fixed such that $f_0=\mathrm{id}$ and tangles $T_2$ and $\tilde T_1=f_1(T_1)$ are congruent up to
dihedral group transformations.\end{dfn}

\centerpic{ball.eps}{A $3$-tangle, its diagram and projection}

Along with the equivalence defined above a weaker equivalence relation is also considered that
allows the legs of the tangle to be moving on the boundary sphere~\cite{Sundberg1998}.

\begin{dfn}\label{weakeq} Tangles $T_1$ and $T_2$ are {\it weakly equivalent\/} if there is a smooth
family $f_t\,{:}\;B^3\to B^3$, $t\in[0,1]$ of smooth homeomorphisms on the ball $B^3$ such that
$f_0=\mathrm{id}$ and $f_1(T_1)=T_2$.\end{dfn}

In \cite{Kanenobu2003} 2-tangles with 7 crossings or less are classified up to weak equivalence.

For graphical representation of tangles conventional {\it diagrams\/} are used which are
non-singular projections onto the equatorial disk $x^2+y^2\le1$, $z=0$ with additional
information at each double point about which ``thread'' is above the other one
(\picref{ball.eps}b). We are more interested here in projections in themselves
(\picref{ball.eps}c).

\begin{dfn} A tangle projection (diagram) is called {\it composite,} if there is a smooth
closed curve inside the boundary circle that intersects the projection exactly twice and encloses
at least one crossing (\picref{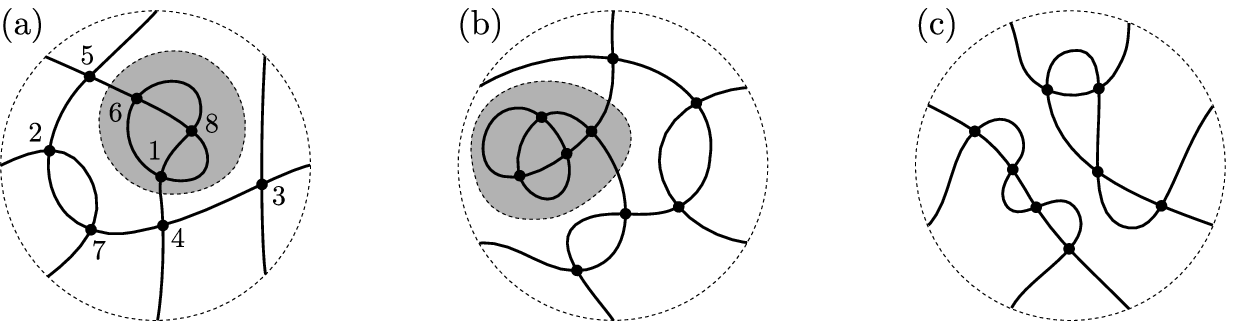}a,b). Otherwise the projection is called {\it
prime.}\end{dfn}

\centerpic{nugatory.eps}{Composite (a,b) and non-connected (c) projections}

As in the case of knots and links we are interested in classification of prime projections only.
Moreover, it is natural to consider only connected projections, i.\,e.~such that it is possible
to get from any crossing to any other moving along the arcs of the projection. An example of
non-connected projection is shown in~\picref{nugatory.eps}c.

We will use the terminology of graph theory for description of projections. In particular, the
crossings of projection will be called {\it vertices,} and arcs will be called {\it edges.} The
set of vertices will be denoted by $V$, and the set of edges will be denoted by $E$.

Finally, for the following we need two more definitions.

\begin{dfn} A crossing (vertex) is called {\it boundary crossing (respectively, vertex)} if some of the
incident edges are connected with the boundary circle.\end{dfn}

\begin{dfn}\label{cut} A {\it cut-crossing (vertex)} is a crossing (vertex) that if removed (together
with the incident edges) produces a non-connected graph.\end{dfn}

For example, the projection \picref{nugatory.eps}a contains five boundary crossings
2,\,3,\,4,\,5,\,7. The crossing number 4 is a cut-crossing.

\section{Cascade diagram of projection}\label{cascade}

Any prime connected projection can be represented as a {\it cascade diagram\/} which is defined
in the following way. Let us construct a system of nested areas in the disk so that there is
exactly one vertex in each annular layer (\picref{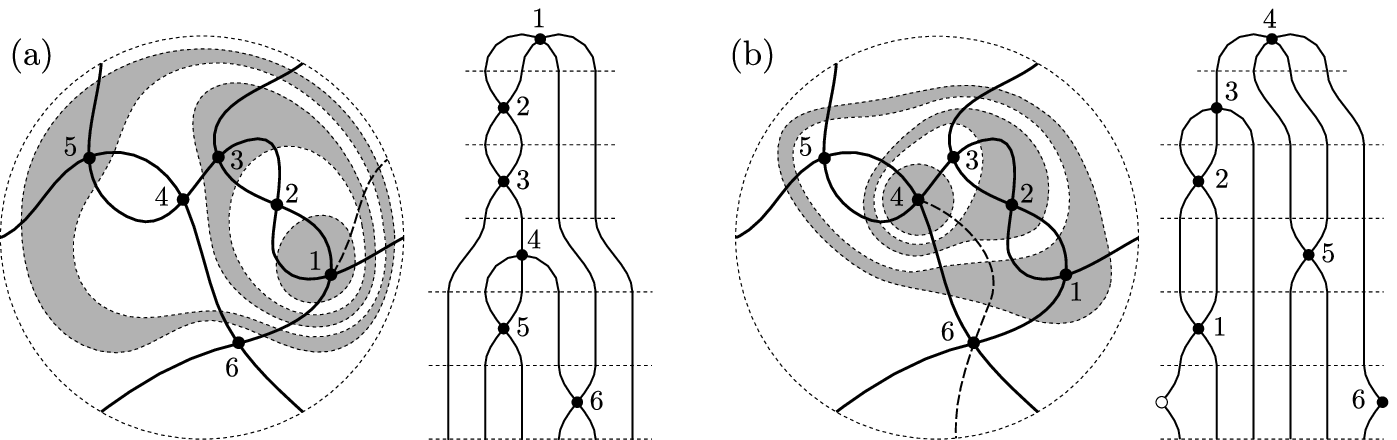}), and a connected sub-projection
is contained inside each area. Of course, it always can be done in many ways.

\centerpic{cascade.eps}{Cascade diagram construction}

Let us cut the disk from the central vertex (vertex~1 in \picref{cascade.eps}) to the boundary
circle so that the cutting line intersects the boundary of each area exactly once (wide-dashed
lines in \picref{cascade.eps}). It is allowed the cut to dissect the projection at vertices only
(see~\picref{cascade.eps}b). Opening the cut, we obtain a representation of the projection as a
falling cascade (\picref{cascade.eps}). We will call this representation {\it cascade diagram} of
the projection.

There is only one vertex at each level of the cascade diagram. Edges incident to this vertex
intersect upper and/or lower boundaries of the level. Only five configurations are possible, that
correspond to the patterns
$$
\OI,\ \O,\ \PI,\ \P,\ \X.
$$
The \OI pattern, which is always the starting pattern, cannot be located at lower levels because
of connectivity of the projection and by construction. Moreover, let us prove the following
statement:

\begin{thm}\label{thm1}
Any prime projection has a cascade diagram free of the \O pattern; any vertex
can be chosen as the starting vertex of such cascade diagram.\end{thm}

We shall prove the theorem by constructing an explicit recurrent algorithm that allows us to
construct a cascade diagram required.

\begin{lmm} Any prime projection has at least two non-cut boundary crossings.\end{lmm}

\centerpic{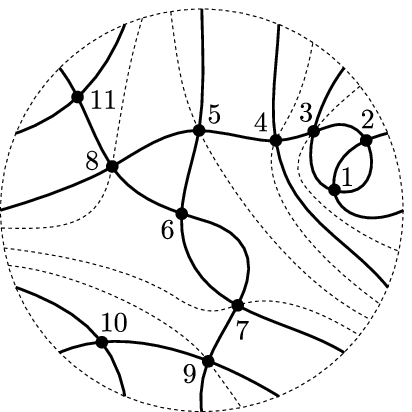}{On the proof of the Lemma}

\begin{proof}
Through each cut-crossing draw a chord (a simple curve with end points on the boundary circle)
that nowhere else intersects the projection (\picref{lemma.eps}). It is always possible by
Def.\,\ref{cut}, and can be done such a way that the chords do not cross each other. Obviously,
there are at least two ``outermost'' chords (e.\,g.~the chords passing through the crossings
3,\,8,\,9 in \picref{lemma.eps}). Any of these chords cuts a ``lune'' (a part of the disk not
containing cut-crossings). But by the cut-crossing definition any lune should contain at least
one crossing; moreover, there is at least one boundary crossing in each lune (contrary is
possible only for composite projections of type~\picref{nugatory.eps}b). So, for each lune there
is at least one non-cut boundary crossing. This completes the proof.\end{proof}

\begin{proof}[Proof of the Theorem~\ref{thm1}]
Set an arbitrary crossing to be the starting crossing of the cascade. We will build the cascade
diagram in bottom-up direction. Choose one of non-cut boundary crossings that is not the starting
one (it is always possible by the Lemma). Put this crossing on the lower $n$-th level
(see~\picref{cascade.eps}). The residual part of the projection is again a connected prime
projection; it has $n{-}1$ crossings. By the Lemma this projection has at least two non-cut
boundary crossings. So again we can choose a crossing different from the starting one and put it
on the $(n{-}1)$-th level. Repeat this procedure until all the crossings (except the starting
one) will arrange; so we construct a cascade diagram required. The theorem is proved.\end{proof}

Thus, any projection can be represented by a cascade diagram containing (besides the starting
pattern \OI) the patterns \PI, \P, \X only; it is easily shown that each of these patterns is
necessary.

\section{Cascade code}

Let us code levels of cascade diagram starting from the top. The pattern \OI is always on the
first level so we will not include it in the code. For other levels we will use the following
coding scheme. Associate a pair $(\alpha_i,m_i)$ with $i$-th level of the cascade diagram. Here
$\alpha_i$ denotes one of the patterns $\P$, $\X$, $\PI$, and an integer $m_i$ determines the
shift of the pattern $\alpha_i$ with respect to the pattern $\alpha_{i-1}$ on the previous level.
We should choose a reference point for each pattern to determine the shift $m_i$: let settle the
origin on the left leg for the $\X$ pattern, and on the central leg for the patterns $\P$ and
$\PI$. The shift $m_i\in0,1,\dots,k_{i-1}{-}1$ counts in counterclockwise direction (in
left-right direction on cascade diagram).

So a diagram with $n$ crossings is coded by an ordered set of $n{-}1$ pairs:
$$
C=\bigl\{(\alpha_2,m_2),\dots,(\alpha_n,m_n)\bigr\}
$$
(recall that the code does not include the starting pattern $\OI$). For example, the code
$$
C=\bigl\{(1,0),(0,5),(0,3),(0,3),(0,5)\bigr\}
$$
corresponds to the cascade diagram \picref{cascade.eps}b.

Cascade code for $k$-tangles with $n$ crossings is almost as compact as the DT-code for knots
with $n$ crossings~\cite{Dowker1983,Rankin2002_1}: we need only 2 bits per crossing additionally
for storing $\alpha_i$. It is properly that the length of cascade code only depends on the
crossing number; it does not depend on the number of tangle components (and therefore on the
number of legs).

Cascade code actually is a recipe how to draw the projection, so we have no problems with
checking of the code drawability (it is well known that it is not an easy task for DT-code). The
last property allows us to use the cascade code effectively for generation and enumeration of
tangle projections with given number of crossings.

Each projection has a lot of cascade representations so the equivalence problem for cascade codes
occurs. In order to resolve this problem it is necessary to define a {\it canonical cascade
code\/} which should correspond uniquely to a given projection. Such code can be constructed in
different ways. We will define it in such a way that the important {\it nesting property\/} will
be satisfied.
\begin{equation}\label{nesting}\raise-\baselineskip\hbox{\vbox{\hsize=0.75\textwidth \it\noindent
If $C_n=\bigl\{(\alpha_2,m_2),\dots,(\alpha_n,m_n)\bigr\}$ is a canonical cascade code, then any
initial segment $C_i=\bigl\{(\alpha_2,m_2),\dots,(\alpha_i,m_i)\bigr\}$, $2\le i\le n{-}1$, of the
code is again the canonical cascade code of the corresponding sub-tangle with $i$~crossings.}}
\tag{$*$}
\end{equation}
In order to do define canonical cascade code, we will at first define an auxiliary invariant code
of projection.

\section{Invariant root-code of projection}

\begin{dfn} A {\it root\/} of a projection is a triple $r=(v,e,f)$ where $v\in V$ is a
vertex of the projection, $e\in E$ is one of four edges incident to $v$, $f$ is one of two faces
incident to $e$ (see examples in~\picref{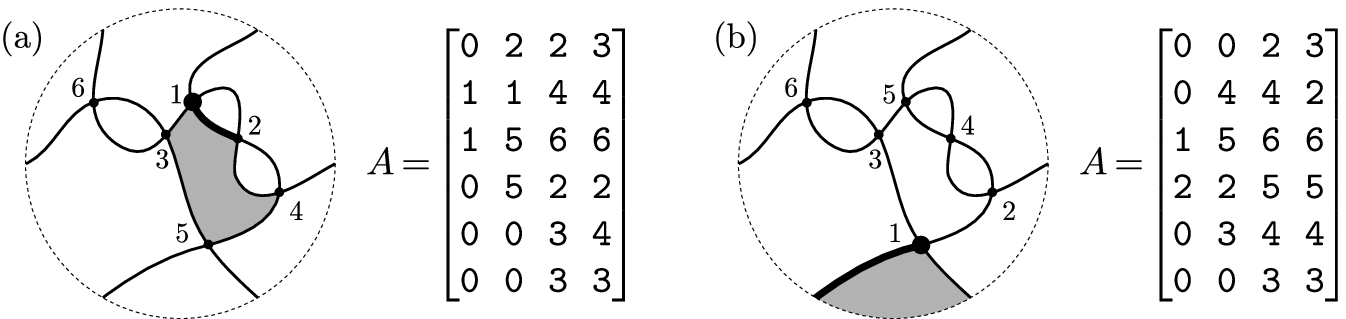}). $v$, $e$, and $f$ is called {\it
root-vertex,} {\it root-edge,} and {\it root-face\/} respectively.\end{dfn}

The location of the root-face relative to the root-edge and the root-vertex determines naturally
a rotation direction (clockwise or counter-clockwise) about the root-vertex which we will call
the {\it labeling direction.}

Vertices of the projection can be uniquely labeled as a root fixed. As an example the following
algorithm may be used:
\begin{quote}\parskip=0mm\ttfamily\footnotesize

0.~~set the number of the root-vertex as 1;

~~~~~~number vertices adjacent to the root-vertex

~~~~~~in the labeling direction starting from the root-edge;

1.~~among the labeled vertices find the vertex $q$

~~~~~~with the smallest number that has unlabeled neighbors;

2.~~assign the numbers to unlabeled vertices around the vertex $q$

~~~~~~in order determined by the labeling direction,

~~~~~~starting from the first unlabeled edge;

3.~~if there are still unlabeled vertices return to the step 1.
\end{quote}

\noindent The vertices of the projection in \picref{rootnum.eps}a,b are labeled using the
algorithm described above.

\centerpic{rootnum.eps}{Two vertex labelings induced by different roots, the adjacency lists;\\
the canonical labeling (b)}

During the labeling algorithm execution an \emph{adjacency list\/} $A_r$ is constructed. The
$i$-th line of the adjacency list contains labels of the vertices adjacent to the $i$-th vertex
listed in the labeling direction. For the boundary vertices we substitute $\mathtt{0}$ at the
positions corresponded to the boundary edges (see~\picref{rootnum.eps}).

The order of numbers in each line of the list is determined up to a cyclic shift. Let us chose
the lexicographically minimal sequence from the four possible variants. (The line
$a=[a_1,\dots,a_n]$ is {\it lexicographically less\/} than line $b=[b_1,\dots,b_n]$ if there is
$j$ such that $a_j<b_j$ and $a_i=b_i$ for any $i<j$.)

Let us agree to write the adjacency list $A_r$ in one line. For example, the root in
\picref{rootnum.eps}a determines the adjacency list
$$
A_r=[\mathtt{0~2~2~3~1~1~4~4~1~5~6~6~0~5~2~2~0~0~3~4~0~0~3~3}].
$$

\begin{dfn} A single-line adjacency list $A_r$ constructed as described above is called the {\it root-code\/}
cor\-re\-spon\-ding to a given root~$r$.\end{dfn}

Thus as a root fixed the corresponding root-code is defined uniquely; in turn a root-code
determines uniquely the corresponding projection (up to rotations and reflections).

We can define a root of a given projection in $8n$ ways. Therefore, there are $8n$ generally
different root-codes of the projection; so in order to define an invariant code we need a rule
that allows to choose a unique root-code from this set. The root-code constructed by that rule we
will call canonical and will denote as $\rcode(\cdot)$.

\begin{dfn} Root-code $\rcode(\cdot)$ is called an {\it invariant root-code\/} if it is
determined uniquely by a given projection:
$$
\rcode(P_1)=\rcode(P_2)\ \Longleftarrow\ P_1\cong P_2.
$$\end{dfn}

For example one of the simplest rules (but in the same time one of the most inefficient) is the
following one: define the invariant root-code as the lexicographic minimum among all $8n$ codes.

In order to construct a canonical cascade code satisfying the nesting property \eqref{nesting} we
need the root-vertex to be {\it boundary\/} and {\it non-cut.} It is obvious that a root-code
will be lexicographically minimal only if the root-vertex has {\it maximum number of legs:} in
this case there will be maximum number of starting zeroes in the root-code. For the projection in
\picref{rootnum.eps}a even these two requirements reduce the search set to two root-vertices 2
and 6.

We can reduce the number of compared codes some more, for example using a face-code notion. Let
us associate an integer with each boundary face --- the number of projection edges adjacent to
the face. Then we can write these numbers for all boundary faces in order defined by the labeling
direction as a string of length~$2k$ (number of legs). The {\it face-code\/} corresponding to a
given (boundary) root is such a sequence that the starting number corresponds to the root-face.
For example the face-code for the root in \picref{rootnum.eps}b is $[\mathtt{2~3~4~4~2~4}]$.

Thus, we select an invariant subset in the set of roots
\begin{equation}\label{rset}
R=\Biggl\{r\;\Biggr|\ \raise-1\baselineskip\vbox{\hsize=17em\noindent
root-vertex is a boundary vertex,\\
root-vertex is a non-cutting vertex,\\
face-code is lexicographically minimal}\Biggr\}.
\end{equation}

Now let us define {\it invariant root-code\/} $\rcode(\cdot)$ as the lexicographically minimal
root-code on the subset $R$:
$$
\rcode(T)=\lexmin_{r\in R}\{A_r(T)\}.
$$

\begin{dfn}\label{canonroot}
The root defining the invariant root-code we will call {\it canonical root.}
\end{dfn}

The subset $R$ for the projection in \picref{rootnum.eps} contains only one root, namely, the
root depicted in \picref{rootnum.eps}b; the invariant root-code is the following:
$$
\rcode(P)=[\mathtt{0~0~2~3~0~4~4~2~1~5~6~6~2~2~5~5~0~3~4~4~0~0~3~3}].
$$

It is important to notice that the average number of elements in the subset $R$ (among all the
projections with given number of crossings) tends to~1 as crossing number $n$ increases. This is
because the fraction of symmetrical face-codes tends to~0. Consequently, in case of large~$n$
with probability close to~1 the canonical root can be found without calculation of the canonical
root-code.

\section{Canonical cascade code}\label{ccc}

Now, using invariant root-code described above, we can define a canonical cascade code satisfying
the nesting property~\eqref{nesting}.

Consider a tangle projection $T_n$ with $n$ vertices. We will construct a canonical cascade
diagram (and canonical cascade code) from bottom to top using a procedure similar to that used in
the proof of the Theorem~\ref{thm1}; let us describe it omit some details.

Find the canonical root-vertex of the projection $P_n$ (construct the canonical root-code if it
is necessary); put this vertex on the $n$-th level of the cascade being constructed. In doing so,
we obtain a tangle projection $P_{n-1}$ (with $n{-}1$ crossings) corresponds to the rest part of
the projection. For the projection $P_{n-1}$ we can again find the canonical root-vertex; and so
the $(n{-}1)$-th level of the cascade diagram will defined. Continuing, we construct a uniquely
determined cascade diagram and, in the same time, a unique sequence of nesting tangle
projections:
\begin{equation}\label{parents}
\OI=P_1\ \lar\ P_2\ \lar\ \cdots\ \lar\ \ P_{n-1} \lar\ P_n.
\end{equation}
The cascade code built using the described algorithm we will call {\it canonical cascade code\/}
of projection.

The canonical cascade for the projection \picref{cascade.eps} is shown in \picref{cascade.eps}b;
the corresponding canonical cascade code is the following:
$$\let\ub=\underbrace
\bigl\{\ub{\ub{\ub{\ub{\ub{(0,0)}_{P_2},(0,0)}_{P_3},(1,1)}_{P_4},(0,5)}_{P_5},(0,2)}_{P_6}\bigr\}.
$$
The sequence \eqref{parents} determines a ``genealogy'' of the projection. For the projection
\picref{cascade.eps} we obtain the series depicted in \picref{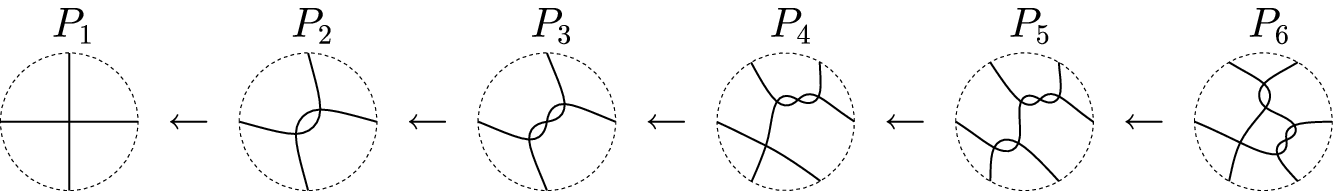}.

\centerpic{geneology.eps}{The genealogy of the projection \picref{cascade.eps}}

\noindent Thus, for a given projection the parent-projection is uniquely defined by the canonical
cascade code. We will use this fact in the next section to construct an enumeration algorithm for
$k$-tangle projections.

\section{Enumeration algorithm}\label{alg}

Algorithms of knots and links classification usually include two stages: 1)~generation, and
2)~sifting. The main objective of the generation stage is to obtain (using some representation of
knot projections, e.\,g.~DT-code) a set that includes all the projections with given number of
crossings. As a rule, the obtained set is redundant, because there are lots of different
representations of the same projection. These duplicates must be found using invariants and
eliminated in the second stage of the algorithm, in the sifting stage. The effectiveness of the
classification algorithm depends on the amount of redundancy and on complexity of the equivalence
checking procedure.

Recent classification algorithms \cite{Hoste2006,Rankin2002_1} realizing the generation stage
based on the succession principle: knots with $n{+}1$ crossings are generating from the (sifted)
set of knots with $n$ crossings. For example, in the
algorithm~\cite{Rankin2002_1,Rankin2002_2,Rankin2002_3} each knot with $n$~crossings generates a
set of child knots with $n{+}1$ crossings using local operations like shown
in~\picref{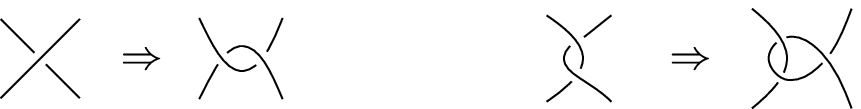}.

\centerpic{local.eps}{Local operations increasing the crossing number by 1~\cite{Rankin2002_1}}

Such approach keeps the amount of duplicates relatively small and prohibits generation of
non-drawable representations that do not correspond to any knot projection (as it take place, as
an example, in the generator based on DT-code).

Obviously, using cascade representation we can construct a generator that satisfies the
properties mentioned above. For a given projection with $2k$ legs the set of child projections is
obtaining by adding one of the symbols $\P$, $\X$, $\PI$ in all possible $2k$ ways to a cascade
diagram of the projection ($6k$ descendants total). The effectiveness of the algorithm can be
enhanced appreciably using canonical cascade code.

Suppose we have the full set $\Pr_n$ (with no repetitions) of projections with $n$ crossings;
where each projection is represented by its canonical cascade code. We can add the
\hbox{$(n{+}1)$-th} vertex (the $(n{+}1)$-th layer) to a given cascade diagram from $\Pr_n$ in
$6k$ ways (where $2k$ is number of legs); so the diagram generates $6k$ descendants with $n{+}1$
crossings. Then the following simple rule can be used to decide if we should reject a given child
projection or tabulate~it.
\begin{equation}\label{rule}\hbox{\vbox{\hsize=0.75\textwidth \it\noindent
If the $(n{+}1)$-th vertex is not canonical root-vertex $\longrightarrow$ reject the projection.}}
\tag{$**$}
\end{equation}
For that first of all it should be checked if the new $(n{+}1)$-th vertex belongs to the
$R$-subset~\eqref{rset} or not. If ``not'', we reject the projection without calculating the
invariant root-code; if ``yes'', and the $R$-subset contains more than one element, the invariant
root-code should be found to identify the canonical root-vertex.

Thus, from the new generation of cascade diagrams only diagrams in canonical form will survive.

The deciding rule~\eqref{rule} ensures the following essential property. Because the canonical
cascade code of a projection determines uniquely its genealogy~\eqref{parents}, projections
descend from different parents are necessarily different. So, in the sifting stage there is no
need to compare projections obtained from different parents.

Moreover, we can get the set $\Pr_{n+1}$ (projections with $n{+}1$ crossings) with no duplicates
at all if we take into account information about symmetries of projections from the $\Pr_n$ set.

First let us show that no asymmetrical projection can generate two equivalent canonical cascade
diagrams.

\begin{thm} Let $P_0$ be a projection with trivial symmetry group; $C_0$ is the canonical
cascade diagram of $P_0$ and cascade diagrams $C_1$ and $C_2$ are descendants of $C_0$
(\picref{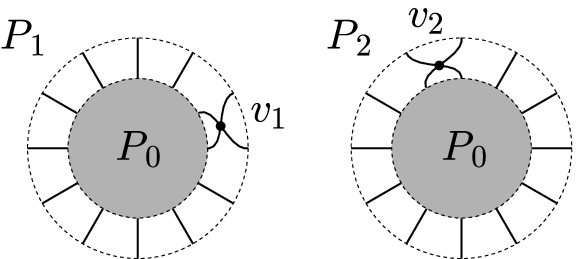}). Let the projections $P_1$ and $P_2$, which correspond to the cascade
diagrams $C_1$ and $C_2$, are equivalent. Then cascade diagrams $C_1$ and $C_2$ cannot be
simultaneously canonical.\end{thm}

\centerpic{proof2.eps}{On the proof}

\begin{proof} Suppose, to the contrary, that cascade diagrams $C_1$ and $C_2$ are both canonical.
Then, by the definition of canonical cascade diagram, both the vertices $v_1$ and $v_2$
(\picref{proof2.eps}) are canonical root-vertices. Now, let us note an obvious corollary of the
canonical root definition~(Def.\,\ref{canonroot}): any transformation from the dihedral group
$D_{2k}$ that transfers a canonical root into a canonical root, transforms the projection into
itself. Therefore, the (non-trivial) transformation $\sigma$ that transfers the root $v_1$ into
the root $v_2$ is in the symmetry group of the projection $P_1\cong P_2$. But then, obviously,
the symmetry group of the projection $P_0$ also contain $\sigma$. This contradiction completes
the proof.
\end{proof}

Thus, we see that only symmetrical projections can produce duplicates; however it is easy to
avoid these duplications. Obviously, two child projections $P_1$ and $P_2$ of a given symmetrical
projection $P_0$ are equivalent if $P_1$ can be transformed into $P_2$ (up to isotopy in the
outside annular layer (see~\picref{proof2.eps})) by a symmetry of $P_0$.
Among all $6k$ child projections we should consider only projections that cannot be transformed
one into another by a transformation mentioned above; let us denote this subset $S$. For example
only two of 12 \P-child projections of projection $P_0$ (\picref{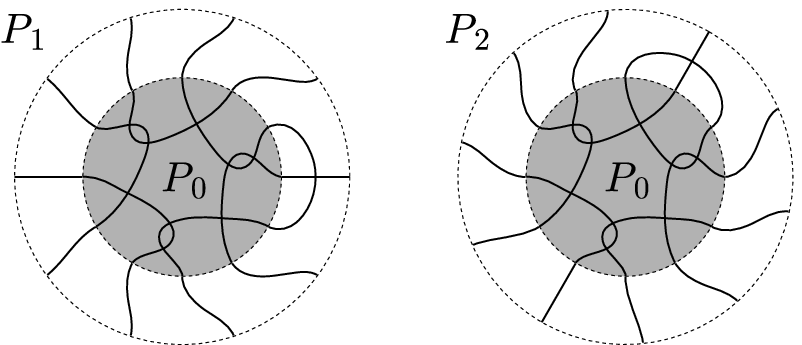}) are in the $S$
subset, namely the projections $P_1$ and $P_2$ in the Figure.

\centerpic{symchild.eps}{Child projections of a symmetrical projection}

Similar to the asymmetrical case, it is possible to proof that if there are equivalent
projections in the subset $S$, then only one of them can be in canonical form.

So, if the rule \eqref{rule} is used (taking symmetries into account) then the set of descendants
of the set $\Pr_n$ is exactly the set $\Pr_{n+1}$.

As evidenced from the above, any projection has a unique parent projection. Therefore, we can
represent the set $\Pr$ of all tangle projections as a genealogical tree starting from the
one-crossing ancestor. In~\picref{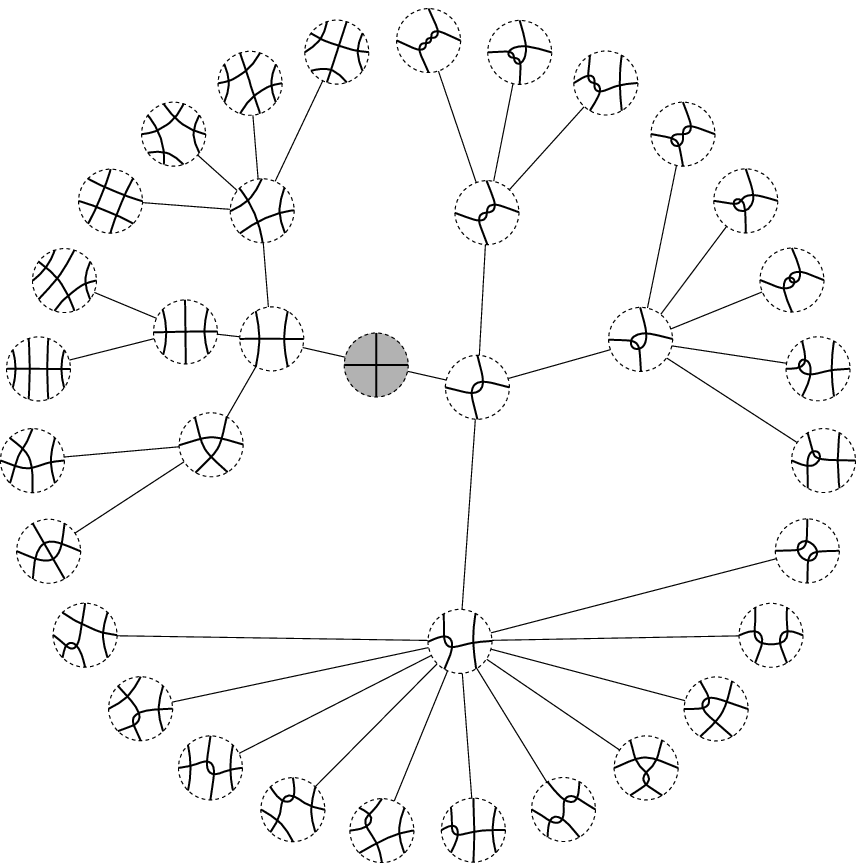} the first three stages of the tree are depicted. It
is easy to demonstrate that the genealogical tree has no dead-end branches: any projection has at
least one descendant.

\centerpic{tree.eps}{Genealogical tree of tangle projections}

In conclusion let us describe the algorithm entirely.

Let $\Pr_n$ be the full set of projections with $n$ crossings. We need to do the following steps
for each projection in $\Pr_n$ to generate the set $\Pr_{n+1}$:
\begin{quote}\parskip=0mm\ttfamily\footnotesize
1.~~find allowable positions for symbols $\P$, $\X$, $\PI$

~~~~~~taking into account symmetry of the projection;

2.~~eliminate the cases when the new $(n{+}1)$-th vertex

~~~~~~is not in the $R$-set~\eqref{rset} of the child projection;

3.~~eliminate the cases when the child projection is composite;

4.~~eliminate the cases when the new $(n{+}1)$-th vertex

~~~~~~is not canonical root-vertex of the child projection;

5.~~add the remaining child projections to the catalogue;

6*.~find flipping-equivalent projections, eliminate the duplicates.
\end{quote}

We will not dwell here on the algorithm of flipping-equivalence checking. The algorithm we have used
differs somewhat from the algorithm~\cite{Rankin2002_1,Rankin2002_2,Rankin2002_3}, we will
describe it in a separate paper.

\section{Tables of tangle projections}

In this section we provide some results of enumeration of tangle projections and some important
sub-classes of the set $\Pr$. We also give tables of alternating $k$-tangles with five crosings
or less.

We use the following notation. As above, $\Pr$ denotes the set of connected prime $k$-tangle
projections. If we assume that projections related by a series of flipes (\picref{flyp.eps}) are
equivalent, we obtain the set of alternating $k$-tangles denoted by $\A$. If we consider
$k$-tangles up to weak equivalence (Def.\,\ref{weakeq}), many projections became equivalent. For
example, any projection that contains a boundary crossing with more than one legs are equivalent
to a projection with smaller crossing number. We call such projections {\it weakly equivalent\/}
and denote the set of weak equivalence classes by $\W$. Finally, we are interested in
classification of $k$-tangle projections that have no more than one edge between any two vertices
(do not contain 2-faces). We call projections of this type {\it reduced\/} and denote the
corresponding set by $\Rr$. Reduced projections play a role analogous to the role of Conway
polyhedra~\cite{Conway1970} in classification of knots and links.

Subsets of projections with fixed crossing number $n$ and number of legs $2k$ we label by
indexes. For example, $\A_{5,4}$ is the set of alternating tangles with 5 crossings and 8 legs;
$\Rr_7$ is the set of all reduced projections with 7 crossings.

Using cascade representation and the algorithm described above we tabulate the sets of
projections $\Pr$, $\A$, $\W$, and $\Rr$ up to 12 crossings.

Table~\ref{table1} contains numbers of projections in subsets $\Pr_{n,k}$ for $n\le12$.

\begin{center}
\texttt{\scriptsize
\begin{tabular}{|c||r|r|r|r|r|r|r|r|r|r|r|r|}
\multicolumn{13}{r}{\normalfont\footnotesize{\bf Table 1 }
Number of tangle projections with $n$ crossings and $2k$ legs
\refstepcounter{table}\label{table1}}\\[1ex]
\hline $k\setminus n$
    & 1 & 2 & 3 & 4 & 5  & 6   & 7      & 8       & 9       & 10         & 11          & 12\\\hline\hline
  2 & 1 & 1 & 2 & 6 & 19 & 71  & 293    & 1\,348  & 6\,568  & 33\,701    & 178\,706    & 973\,085\\
  3 & . & 1 & 2 & 8 & 29 & 138 & 638    & 3\,237  & 16\,805 & 90\,239    & 494\,151    & 2\,756\,453\\
  4 & . & . & 2 & 8 & 41 & 210 & 1\,125 & 6\,138  & 34\,112 & 192\,278   & 1\,096\,560 & 6\,317\,363\\
  5 & . & . & . & 5 & 31 & 231 & 1\,458 & 9\,183  & 56\,084 & 340\,885   & 2\,060\,224 & 12\,446\,400\\
  6 & . & . & . & . & 16 & 161 & 1\,406 & 10\,572 & 74\,331 & 499\,902   & 3\,276\,104 & 21\,112\,641\\
  7 & . & . & . & . & .  & 60  & 840    & 8\,818  & 75\,747 & 591\,091   & 4\,327\,816 & 30\,451\,898\\
  8 & . & . & . & . & .  & .   & 261    & 4\,702  & 56\,199 & 541\,570   & 4\,628\,641 & 36\,633\,417\\
  9 & . & . & . & . & .  & .   & .      & 1\,243  & 26\,753 & 361\,106   & 3\,846\,580 & 35\,758\,786\\
  10& . & . & . & . & .  & .   & .      & .       & 6\,257  & 155\,593   & 2\,332\,512 & 27\,199\,662\\
  11& . & . & . & . & .  & .   & .      & .       & .       & 32\,721    & 916\,595    & 15\,123\,600\\
  12& . & . & . & . & .  & .   & .      & .       & .       & .          & 175\,760    & 5\,464\,661\\
  13& . & . & . & . & .  & .   & .      & .       & .       & .          & .           & 963\,900\\\hline
all & 1 & 2 & 6 & 27& 136& 871 & 6\,021 & 45\,241 & 352\,856& 2\,839\,086& 23\,333\,649& 195\,201\,866\\\hline
\end{tabular}}
\end{center}

For the set of alternating $k$-tangles along with the enumeration results (Table~\ref{table2}) we
provide also tables of tangle pictures \picref{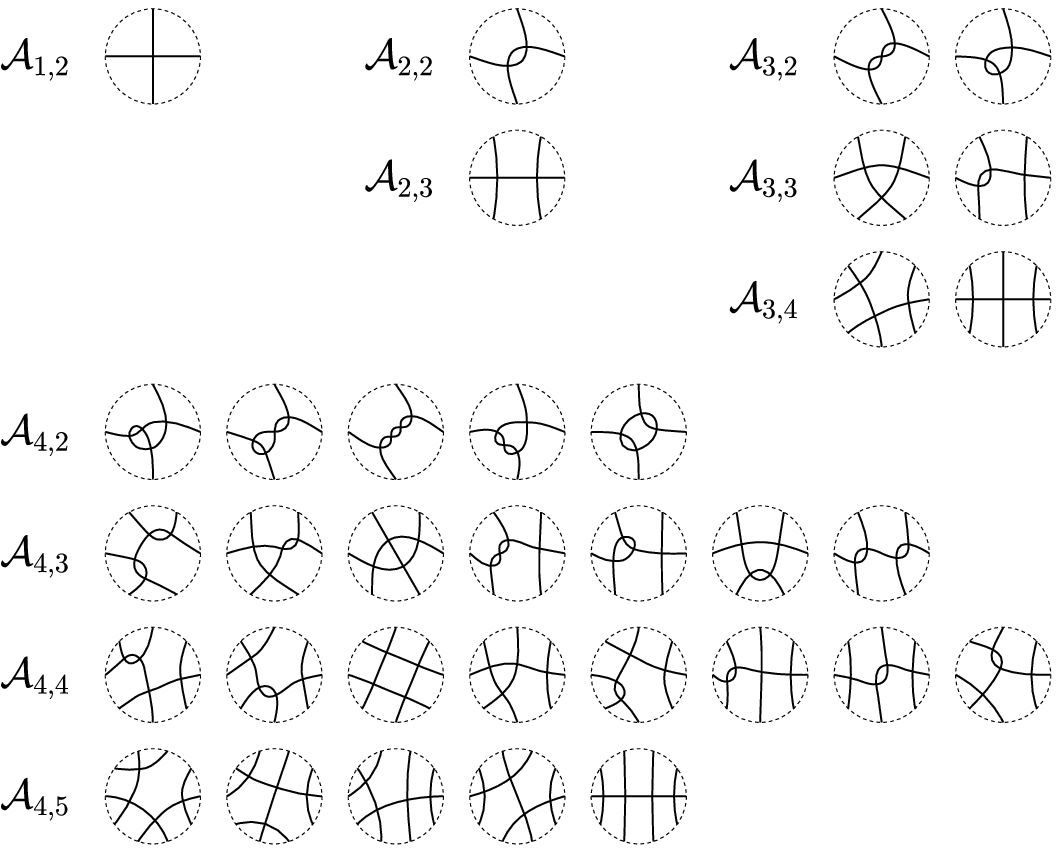} and \picref{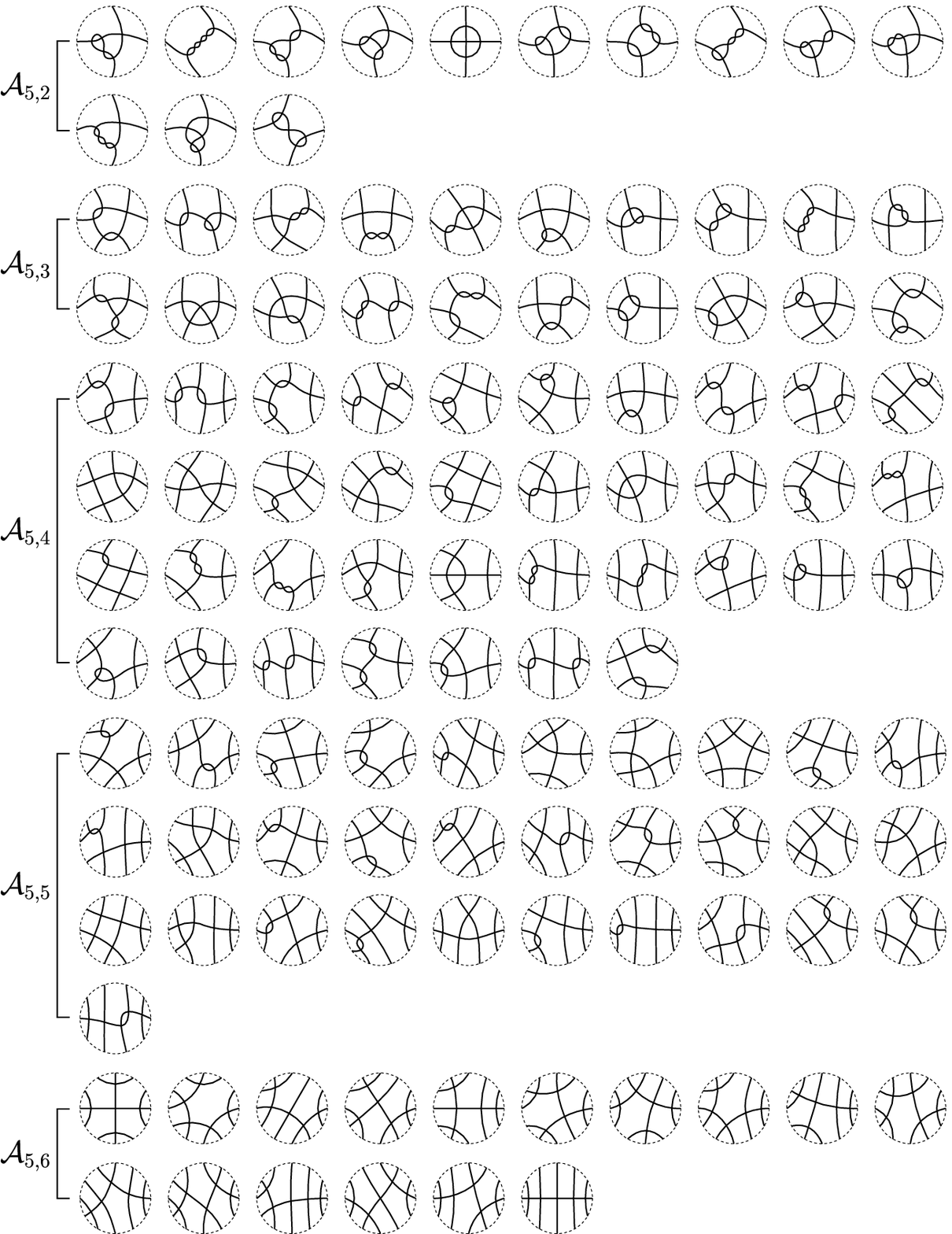}.

\centerpic{tangles14.eps}{Alternating tangles with up to 4 crossings}

\centerpic{tangles5.eps}{Alternating tangles with 5 crossings}

\begin{center}
\texttt{\scriptsize
\begin{tabular}{|c||r|r|r|r|r|r|r|r|r|r|r|r|}
\multicolumn{13}{r}{\normalfont\footnotesize{\bf Table 2 }
Number of alternating tangles with $n$ crossings and $2k$ legs
\refstepcounter{table}\label{table2}}\\[1ex]
\hline $k\setminus n$
    & 1 & 2 & 3 &  4 &   5 &   6 &      7 &       8 &        9 &          10 &           11 &            12 \\
\hline\hline
2   & 1 & 1 & 2 &  5 &  13 &  36 &    111 &     373 &   1\,362 &      5\,378 &      22\,807 &      102\,617 \\
3   & . & 1 & 2 &  7 &  20 &  77 &    276 &  1\,135 &   4\,823 &     21\,734 &     101\,307 &      488\,093 \\
4   & . & . & 2 &  8 &  37 & 157 &    687 &  3\,052 &  13\,981 &     65\,797 &     317\,506 &   1\,565\,163 \\
5   & . & . & . &  5 &  31 & 209 & 1\,128 &  5\,986 &  30\,556 &    155\,964 &     795\,918 &   4\,092\,027 \\
6   & . & . & . &  . &  16 & 161 & 1\,294 &  8\,528 &  51\,475 &    294\,366 &  1\,637\,855 &   8\,979\,493 \\
7   & . & . & . &  . &   . &  60 &    840 &  8\,206 &  62\,895 &    428\,254 &  2\,702\,902 &  16\,313\,106 \\
8   & . & . & . &  . &   . &   . &    261 &  4\,702 &  52\,815 &    460\,189 &  3\,475\,551 &  23\,979\,733 \\
9   & . & . & . &  . &   . &   . &      . &  1\,243 &  26\,753 &    341\,878 &  3\,327\,424 &  27\,625\,056 \\
10  & . & . & . &  . &   . &   . &      . &       . &   6\,257 &    155\,593 &  2\,221\,544 &  23\,869\,621 \\
11  & . & . & . &  . &   . &   . &      . &       . &        . &     32\,721 &     916\,595 &  14\,473\,275 \\
12  & . & . & . &  . &   . &   . &      . &       . &        . &           . &     175\,760 &   5\,464\,661 \\
13  & . & . & . &  . &   . &   . &      . &       . &        . &           . &            . &      963\,900 \\\hline
all & 1 & 2 & 6 & 25 & 117 & 700 & 4\,597 & 33\,225 & 250\,917 & 1\,961\,874 & 15\,695\,169 & 127\,916\,745 \\\hline
\end{tabular}}
\end{center}

\newpage

For the subset $\W$ of weak equivalence classes of projections we obtain the following results.

\centerpic{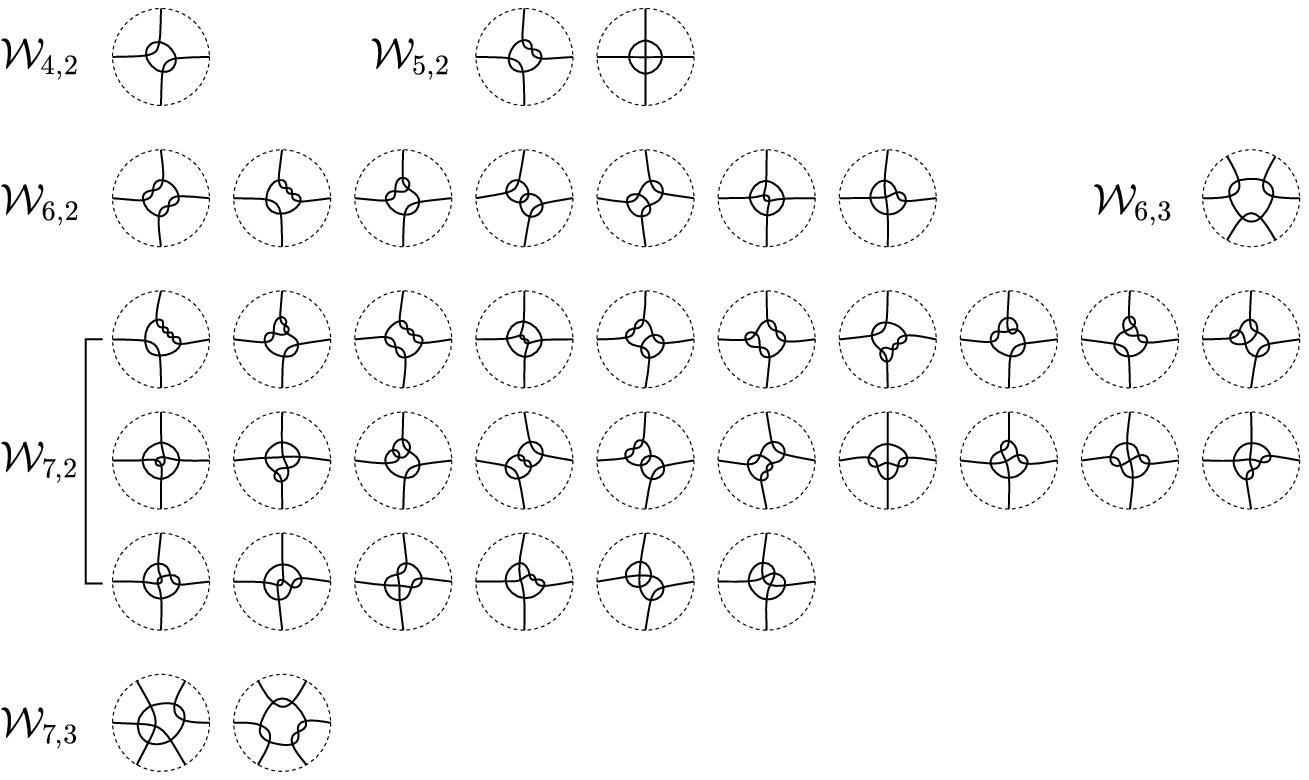}{Alternating tangles with 7 crossings or less up to weak equivalence}

Note, that this catalog agrees with the table presented in~\cite{Kanenobu2003}.

Finally, let us present the results for the subset $\Rr$ of reducing projections.

\begin{center}
\texttt{\scriptsize
\begin{tabular}{|c||r|r|r|r|r|r|r|r|r|r|r|r|}
\multicolumn{13}{r}{\normalfont\footnotesize{\bf Table 4 }
Number of reduced projections with $n$ crossings and $2k$ legs
\refstepcounter{table}\label{table4}}\\[1ex]
\hline $k\setminus n$
    & 1 & 2 & 3 & 4 & 5  &  6 &   7 &     8 &      9 &      10 &      11 &         12 \\\hline\hline
  2 & 1 & . & . & . &  1 &  1 &   3 &     9 &     26 &      74 &     238 &        770 \\
  3 & . & 1 & 1 & 1 &  1 &  4 &   7 &    24 &     69 &     226 &     719 &     2\,423 \\
  4 & . & . & 2 & 2 &  4 &  7 &  21 &    58 &    185 &     596 &  1\,998 &     6\,753 \\
  5 & . & . & . & 5 &  9 & 22 &  49 &   152 &    458 &  1\,545 &  5\,188 &    17\,990 \\
  6 & . & . & . & . & 16 & 42 & 126 &   355 & 1\,144 &  3\,769 & 13\,012 &    45\,515 \\
  7 & . & . & . & . & .  & 60 & 228 &   799 & 2\,586 &  8\,850 & 30\,754 &   109\,843 \\
  8 & . & . & . & . & .  & .  & 261 &1\,288 & 5\,164 & 18\,745 & 68\,142 &   248\,891 \\
  9 & . & . & . & . & .  & .  & .   &1\,243 & 7\,525 & 33\,856 &134\,834 &   520\,884 \\
  10& . & . & . & . & .  & .  & .   & .     & 6\,257 & 44\,482 &222\,482 &   962\,620 \\
  11& . & . & . & . & .  & .  & .   & .     & .      & 32\,721 &266\,270 &1\,464\,500 \\
  12& . & . & . & . & .  & .  & .   & .     & .      & .       &175\,760 &1\,607\,405 \\
  13& . & . & . & . & .  & .  & .   & .     & .      & .       & .       &   963\,900 \\\hline
all & 1 & 1 & 3 & 8 & 31 &136 & 695 &3\,928 &23\,414 &144\,864 &919\,397 &5\,951\,494 \\\hline
\end{tabular}}
\end{center}

\end{document}